\documentclass[12pt]{amsart}
\usepackage{dsfont} 
\usepackage{exscale} 

\newtheorem{thm}{Theorem}[section]
\newtheorem{cor}[thm]{Corollary}
\newtheorem{lem}[thm]{Lemma}
\newtheorem{prop}[thm]{Proposition}
\theoremstyle{definition}
\newtheorem{defn}[thm]{Definition}
\theoremstyle{remark}
\newtheorem{rem}[thm]{Remark}
\numberwithin{equation}{section}
\DeclareMathOperator{\DET}{det} 
\DeclareMathOperator{\CHAR}{char} 
\DeclareMathOperator{\KER}{Ker} 
\DeclareMathOperator{\SPEC}{Spec} 
\DeclareMathOperator{\GAL}{G} 
\DeclareMathOperator{\CODIM}{codim} 
\DeclareMathOperator{\HF}{H} 
\DeclareMathOperator{\AUT}{Aut} 
\DeclareMathOperator{\cd}{cd} 

\providecommand{\llbracket}{[\![} 
\providecommand{\rrbracket}{]\!]} 

\newcommand*{\mlsymbol}[1]{\mathit{#1}} 
\newcommand{\ignore}[1]{} 

\newcommand*{\conjunct}{\circ} 
\newcommand*{\isomorph}{\cong} 
\newcommand*{\LeftD}{\mathbf{L}} 
\newcommand*{\RightD}{\mathbf{R}} 
\newcommand*{\To}{\longrightarrow} 

\newcommand*{\Complex}{\mathbb C} 
\newcommand*{\Rat}{\mathbb Q} 
\newcommand*{\Int}{\mathbb Z} 
\newcommand*{\Zp}[1][p]{\Int_{#1}} 
\newcommand*{\Cp}[1][p]{\Complex_{#1}} 

\newcommand*{\ideal}[1]{\mathfrak{#1}} 

\newcommand*{\set}[1]{\left\{ #1 \right\}}

\newcommand*{\card}[1]{\sharp#1} 
\newcommand*{\ggt}[2]{({#1},{#2})} 

\newcommand*{\Dsum}{\bigoplus} 
\newcommand*{\tensorp}[1]{\otimes_{#1}} 
\newcommand*{\dtensorp}[1]{\tensorp{#1}^{\mathbb{L}}} 
\newcommand*{\tensor}{\otimes} 
\newcommand*{\algebra}[2]{{#1}[ #2 ]} 
\newcommand*{\funcfield}[2]{#1 ( #2 )} 
\newcommand*{\power}[2]{#1 \llbracket #2 \rrbracket} 
\newcommand*{\cmplx}[1]{#1^{\bullet}} 
\newcommand*{\mult}[1]{#1^\times} 
\newcommand*{\algclosure}[1]{\overline{#1}} 

\newcommand*{\quotient}[1]{Q(#1)} 
\newcommand*{\cart}{\mathfrak{I}} 
\newcommand*{\Pic}{\mlsymbol{Pic}} 
\newcommand*{\Gm}{\mathbb{G}_m} 

\newcommand*{\Tor}[2]{\operatorname{Tor}_{#2}^{#1}} 
\newcommand*{\Hc}{\HF} 
\newcommand*{\etH}{\HF_{\mlsymbol{\acute{e}t}}} 
\newcommand*{\IwaH}{\HF_{\mlsymbol{Iw}}} 
\newcommand*{\etSect}{\Gamma_{\mlsymbol{\acute{e}t}}} 
\newcommand*{\IwaSect}{\Gamma_{\mlsymbol{Iw}}} 

\newcommand*{\fG}{\pi_1^{\mlsymbol{\acute{e}t}}} 
\newcommand*{\Order}[1][]{\mathcal{O}_{#1}} 
\newcommand*{\ValR}[1][p]{\mathcal{O}_{#1}} 
\newcommand*{\CyclF}[1]{\funcfield{\Rat}{\pru{#1}}} 
\newcommand*{\Rc}{\RightD\IwaSect/c} 

\newcommand*{\pru}[1]{\zeta_{#1}} 
\newcommand*{\cyclchar}{\varepsilon_{\mlsymbol{cycl}}} 
\newcommand*{\geoF}{\mathcal{F}} 

\newcommand*{\id}{\mlsymbol{id}} 
\newcommand*{\pr}[1]{p_{#1}} 
\newcommand*{\Tw}[1]{\mlsymbol{T\!w}_{#1}} 
\newcommand*{\rec}{\mlsymbol{rec}} 

\newcommand*{\prKL}[1][]{\psi_{#1}} 
\newcommand*{\Ind}[1][K_{\infty}/\Rat]{\operatorname{Ind}_{#1}} 
\newcommand*{\Stickel}[1]{\xi_{#1}} 

\begin{document}

\title[Equivariant main conjecture]{On the Equivariant Main Conjecture of Iwasawa Theory}
\author{Malte Witte}%
\address
{%
Mathematisches Institut\\
Universit\"at Leipzig\\
Postfach 920\\
D-04109 Leipzig\\
Germany
}%
\email{witte@mathematik.uni-leipzig.de}%


\date{\today}%
\begin{abstract}
Refining results of David Burns and Cornelius Greither,
respectively Annette Huber and Guido Kings, we formulate and prove
an equivariant version of the main conjecture of Iwasawa theory
for abelian number fields.
\end{abstract}
\maketitle

\section*{Introduction}

The main conjecture of Iwasawa theory for an abelian number field
in its classical formulation describes the Galois-module structure
of the class groups in the limit over the intermediate fields of
its cyclotomic $\Zp$-extension. The eigenspace of this limit with
respect to a Dirichlet character $\chi$ associated to the ground
field is related to the corresponding $p$-adic $L$-function or to
the eigenspace of the group of global units modulo cyclotomic
units, depending on the parity of $\chi$. The conjecture was
proved in 1986 by B.~Mazur and A.~Wiles \cite{MW1}.

Recently, D.~Burns and C.~Greither \cite{BG1} deduced an
equivariant version of the main conjecture as the key to their
proof of the equivariant Tamagawa number conjecture. Here,
`equivariant' refers to the fact that one retains the full Galois
module structure instead of decomposing the modules by characters.

A.~Huber and G.~Kings \cite{HK1} also use a variant of the Iwasawa
main conjecture in their proof of the Tamagawa number conjecture
for Dirichlet motives. It consists, like the classical
formulation, of a separate statement for each Dirichlet character.
In particular, it is weaker than the formulation in \cite{BG1}.

In the present article, we use this statement and the Theorem of
Ferrero-Washington to reprove the equivariant conjecture of
\cite{BG1} in a slightly more general form.

Fix an odd prime $p$ and let $(K_{\infty},\rho,U)$ be a triple
consisting of the cyclotomic extension $K_\infty$ of an abelian
number field, a one-dimensional representation $\rho$ of
$\GAL(\algclosure{\Rat}/\Rat)$ on a finite extension $\ValR$ of
$\Zp$, and an open subscheme $U$ of $\SPEC \Int$, subject to the
condition that the ramification index in $K_\infty/\Rat$ of every
point in $U$ is prime to $p$.

Like \cite{BG1} and \cite{HK1} we use continuous \'etale
cohomology to assign to each of such triples a complex
$\cmplx{R}=\RightD\IwaSect(U_{K_\infty},\ValR(\rho))$ of modules
over the profinite group ring
$\Omega=\power{\ValR}{\GAL(K_\infty/\Rat)}$. Further, we define a
cyclotomic element $c(U_{K_\infty},\rho)$ in the first cohomology
module of $\cmplx{R}$ and a $p$-adic $L$-element
$L(U_{K_\infty},\rho^{-1}\cyclchar)$ in the quotient ring of
$\Omega$. Here, $\cyclchar$ denotes the cyclotomic character. The
quotient $\cmplx{R}/\Omega c(U_{K_\infty},\rho)$ turns out to be a
perfect complex that is torsion, i.e. acyclic after base change to
the quotient ring. Note that this is no longer true if we drop the
condition on the type of ramification in $U$.

Using the determinant functor of F.~Knudsen and D.~Mumford we can
attach to each perfect torsion complex $\cmplx{P}$ an invertible
fractional ideal of $\Omega$ called the characteristic ideal of
$\cmplx{P}$. Our main result then reads as follows.

\begin{thm}[see Theorem \ref{thm:Main Theorem}]\label{thm:premain}\
\begin{description}
\item[(i) Vanishing of the $\mu$-invariant] Let $\ideal{p}$ be a
prime ideal of codimension $1$ of $\Omega$, with $p\in\ideal{p}$.
Then $(\cmplx{R}/\Omega c(U_{K_\infty},\rho))_{\ideal{p}}$ is
acyclic.

\item[(ii) Iwasawa main conjecture]
$L(U_{K_\infty},\rho^{-1}\cyclchar)$ generates the characteristic
ideal of $\cmplx{R}/\Omega c(U_{K_\infty},\rho)$.
\end{description}
\end{thm}

The formulation of the main conjecture in \cite{BG1} corresponds
to Theorem~\ref{thm:premain}.(ii) for all triples
$(\CyclF{np^{\infty}},\cyclchar^r,\SPEC \algebra{\Int}{1/np})$,
the version in \cite{HK1} to triples
$(\Rat_\infty,\chi\cyclchar^r,\SPEC \algebra{\Int}{1/p})$. Here,
$r$ and $n$ are integers, $\Rat_\infty$ is the $\Zp$-extension of
$\Rat$ and $\chi$ is any Dirichlet character.

The relation to the classical Iwasawa main conjecture is
established by the fact that the first and second cohomology
modules of the complex $\cmplx{R}/\Omega c(U_{K_\infty},\rho)$ for
$\rho=\cyclchar$ is essentially given by the limit  of the
$p$-primary parts of the global units modulo cyclotomic units,
respectively of the class groups, taken over the intermediate
fields of $K_{\infty}/\Rat$. Using
$\cmplx{\tilde{R}}/c(U_{K_\infty},\rho)$ in lieu of these
classical objects leads to a smoother formulation of the
conjecture that circumvents the problems usually connected to $p$
dividing the order of $\GAL(K/\Rat)$ (see the discussion in
\cite{HK1}).

The main idea of the proof of Theorem~\ref{thm:premain} is
essentially the same as in \cite{BG1}. However, we can clarify the
argument considerably by using the result of \cite{HK1}.
Originally, D.~Burns and C.~Greither derived their theorem from
the result of B.~Mazur and A.~Wiles. This approach necessitates
some rather involved deduction steps to deal with the first
cohomology group of $\cmplx{R}/\Omega c(U_{K_\infty},\rho)$, in
particular for Theorem~\ref{thm:premain}.(i). (Note that this step
of the argument in \cite{BG1} contains an inaccuracy that was
later corrected in the appendix of \cite{Fl1}.) The additional
strength of the main conjecture in \cite{HK1} allows us to present
a comparatively quick proof of this part of the theorem. Recall
that A.~Huber and G.~Kings do not use the result of \cite{MW1}.
Instead, they give an independent proof of their statement, using
the Euler system approach of V.~A.~Kolyvagin and K. Rubin
\cite{Ru1}.

As \cite{BG1} and \cite{HK1}, we do not treat the case $p=2$, but
this gap has meanwhile been filled by M.~Flach in \cite{Fl1}.

The article is organised as follows. In
Section~\ref{sec:CharIdeals} we introduce the characteristic ideal
of a perfect torsion complex. Section~\ref{sec:ProfGroupRing}
consists of a collection of algebraic properties of $\Omega$ that
turn out to be useful in the later sections.

The definition of the complex
$\RightD\IwaSect(U_{K_\infty},\ValR(\rho))$ is given in
Section~\ref{sec:IwasawaCohomology}. In the subsequent section we
calculate its cohomology modules in the special case
$\rho=\cyclchar$, which is closely related to classical Iwasawa
theory.

To deal with the ramification of $\rho$ we need an explicit
description of the relative cohomology modules associated to
closed subschemes of $U$. This is achieved in
Section~\ref{sec:LocalFactors}.

In Section~\ref{sec:CyclotomicElements} we extend the classical
construction of cyclotomic elements and $L$-elements to our
setting. The final section is devoted to the proof of the main
theorem.

\subsection*{Acknowledgements}

The paper was partially written up during a visit at the Max
Planck Institute for Mathematics in Bonn. The author would like to
thank this institution for its hospitality and support. Before all
others, he wishes to express his gratitude to A.~Huber for
introducing him to this subject and for numerous valuable
discussions.

\section{Characteristic Ideals}\label{sec:CharIdeals}

The notion of the characteristic ideal of a perfect torsion
complex is a variant of the usual determinant functor of
F.~Knudsen and D.~Mumford \cite{KM1}. It is less flexible than the
latter, but easier to handle.

Let $R$ be any commutative ring (with unit) and denote by
$\quotient{R}$ the total ring of fractions of $R$. Further, we
write $\cart(R)$ for the abelian group of invertible fractional
ideals, i.e. $R$-submodules $I$ of $\quotient{R}$ which are
locally free of rank $1$ and which satisfy
$\quotient{R}\tensorp{R}I=\quotient{R}$.

We can view $\cart$ as a functor from the category of commutative
rings to abelian groups if we restrict the morphisms of the former
to the following class.
\begin{defn}
We call a ring homomorphism $\phi: R\To S$ extendable if it
extends to a homomorphism $\quotient{R}\To\quotient{S}$, also
denoted by $\phi$.
\end{defn}
Examples of extendable homomorphisms include all flat
homomorphisms and all integral extensions.

If $\phi: R\To S$  is extendable, then $\cart(\phi)$ is given by
$$
\cart(\phi)(I)=\phi(I)S
$$
for all $I\in\cart(R)$.

Assume that $R$ is noetherian. Then an element of $\cart(R)$ is
uniquely determined by the following local conditions.

\begin{prop}\label{prop:comparison of ideals}
Let $R$ be noetherian and $I,J\in\cart(R)$. Then $I=J$ if and only
if $IR_\ideal{p}=JR_\ideal{p}$ for all nonzerodivisors $r$ and all
primes $\ideal{p}$ associated to $R/rR$. These are exactly the
primes of codimension $1$ if $R$ is Cohen-Macaulay.
\end{prop}
\begin{proof}
This follows by the same argument as \cite{Ei1}, Prop. 11.3.
\end{proof}

\begin{defn}
We call a complex $\cmplx{P}$ of $R$-modules a torsion complex if
$\quotient{R}\tensorp{R}\cmplx{P}$ is acyclic. $\cmplx{P}$ is
called perfect if it is quasi-isomorphic to a bounded complex of
finitely generated projective $R$-modules.
\end{defn}

Let $\DET_R \cmplx{P}$ denote the determinant of $\cmplx{P}$
according to F.~Knudsen and D.~Mumford \cite{KM1}. If $\cmplx{P}$
is a perfect torsion complex, then the natural isomorphism
$$\quotient{R}\tensorp{R}\DET_R\cmplx{P}=\quotient{R}$$
allows us to view $\DET_R\cmplx{P}$ as an invertible fractional
ideal of $R$.

\begin{defn}
We call $\CHAR \cmplx{P} =(\DET_R\cmplx{P})^{-1}\in\cart(R)$ the
characteristic ideal of $\cmplx{P}$.
\end{defn}

The characteristic ideal enjoys the following properties.

\begin{prop}\label{prop:properties of char}
Let $\cmplx{P}$ be a perfect torsion complex of $R$-modules.
\begin{enumerate}
\item $\CHAR \cmplx{P}$ depends only on the quasi-isomorphism
class of $\cmplx{P}$.

\item $\CHAR \cmplx{P}[1]=(\CHAR \cmplx P)^{-1}$.

\item If $\cmplx{P}_1\To\cmplx{P}_2\To\cmplx{P}_3$ is a
distinguished triangle of perfect torsion complexes in the derived
category, then
$$\CHAR \cmplx{P}_2=\CHAR\cmplx{P}_1\CHAR\cmplx{P}_3.$$

\item If $\phi:R\To S$ is an extendable homomorphism, then
$\LeftD\phi_*(\cmplx{P})=S\dtensorp{R}\cmplx{P}$ is a perfect
torsion complex of $S$-modules and
$$\CHAR\LeftD\phi_*(\cmplx{P})=\cart(\phi)(\CHAR\cmplx{P}).$$

\item If the cohomology modules of $\cmplx{P}$ are themselves
perfect, i.e. of finite $\Tor{}{}$-dimension, then
$$\CHAR\cmplx{P}=\prod_{n\in\Int}(\CHAR \HF^n\cmplx{P})^{(-1)^n}.$$

\item If $R$ is a noetherian and normal domain and $M$ any torsion
module of finite projective dimension (considered as complex
concentrated in degree $0$), then $\CHAR M$ coincides with the
content of $M$, as defined in \cite{Bu1}, VII, \S4.5. In
particular, if $R=\power{\Zp}{T}$, then $\CHAR M$ is the
characteristic ideal of Iwasawa theory.
\end{enumerate}
\end{prop}
\begin{proof}
Everything follows easily from the corresponding properties of the
determinant functor, as given in \cite{KM1}.
\end{proof}

\begin{rem}
If $R$ is not reduced, then the usual determinant functor is
additive only for the class of `true' triangles. In the following,
we will only consider reduced rings. However, note that in our
setting, (iii) is indeed true for arbitrary distinguished
triangles. The reason is that one can always replace the
distinguished triangle by a true triangle of strictly perfect
torsion complexes, the particular choice of which, according to
(i), does not matter. For the determinant functor, it is this
non-canonical choice that causes trouble.
\end{rem}

\begin{rem}
One can also deduce Proposition \ref{prop:properties of char} from
the results of \cite{BreuBurns} on the more sophisticated notion
of the refined Euler characteristic. To this end, note that for a
perfect torsion complex $\cmplx{P}$, the only trivialisation is
the zero map and $\CHAR \cmplx{P}$ is the image of
$-\chi(\cmplx{P},0)$ under the natural homomorphism
$K_0(R,\quotient{R})\To\cart(R)$.
\end{rem}

\section{The Profinite Group Ring of a
$\Zp$-Extension}\label{sec:ProfGroupRing}

In this section we will assemble some useful facts about
cyclotomic $\Zp$-extensions and profinite group rings. A large
part of the material can also be found in \cite{BG1}, \S6.1.

Throughout this article, $p$ will denote a fixed odd prime. Let
$\Rat_\infty$ be the unique $\Zp$-extension of $\Rat$. The
cyclotomic $\Zp$-extension of a number field is given by
$K_\infty=K\Rat_\infty$. We shall always make the additional
assumption that $K$ is an abelian extension of $\Rat$. The Theorem
of Kronecker-Weber then shows that there exists the following
distinguished choice of subfields of $K_\infty$.

\begin{defn}
Let $K_0\subset K_\infty$ be the subfield that is uniquely
determined by the following two properties.
\begin{enumerate}
\item
$\GAL(K_\infty/\Rat)=\GAL(K_0/\Rat)\times\GAL(\Rat_\infty/\Rat)$,

\item $p^2$ does not divide the conductor of $K_0$.
\end{enumerate}
We write $K_n$ for the subfield of $K_\infty$ of degree $p^n$ over
$K_0$.
\end{defn}

Let $\ValR$ be the valuation ring of an arbitrary finite extension
of $\Rat_p$ and write
$$
\Omega=\power{\ValR}{\GAL(K_\infty/\Rat)}=\varprojlim_n\algebra{\ValR}{\GAL(K_n/\Rat)}
$$
for the profinite group ring with coefficients in $\ValR$. Assume
for simplicity that $\ValR$ contains all values of the characters
of $\GAL(K_0/\Rat)$. If $P_\infty$ is the maximal $p$-extension of
$\Rat$ inside $K_\infty$, then
$$
\Omega\isomorph\prod_{\theta}\power{\ValR}{\GAL(P_\infty/\Rat)},
$$
where the product runs through the characters $\theta$ of
$\GAL(K_\infty/P_\infty)$. Observe that
$\power{\ValR}{\GAL(P_\infty/\Rat)}$ is a local Cohen-Macaulay
ring of Krull dimension $2$, but it is not regular unless
$P_\infty=\Rat_\infty$.

The normalisation of $\Omega$ in its total quotient ring
$\quotient{\Omega}$ is given by
$$
\widetilde{\Omega}\isomorph\prod_{\chi}\power{\ValR}{\GAL(\Rat_\infty/\Rat)}.
$$
Here, $\chi$ runs through the characters of $\GAL(K_0/\Rat)$. Note
that
$$
\algebra{\Omega}{1/p}=\algebra{\widetilde{\Omega}}{1/p}.
$$

The prime ideals $\ideal{p}$ of codimension $1$ of $\Omega$ with
$p\in\ideal{p}$ play a special role in our considerations. Recall
that a torsion module $M$ over $\power{\Zp}{T}$ has vanishing
Iwasawa $\mu$-invariant if $M$ is finitely generated as
$\Zp$-module. We generalise this as follows.

\begin{lem}\label{lem:criterion for vanishing}
Let $M$ be an $\Omega$-module which is finitely generated as
$\ValR$-module. Then $M_{\ideal{p}}=0$ for all prime ideals
$\ideal{p}$ of codimension $1$ containing $p$.
\end{lem}
\begin{proof}
We can view $M$ as a module over
$$
\power{\Zp}{T}\isomorph\power{\Zp}{\GAL(K_\infty/K_0)}
$$
via the natural inclusion
$$
i:\power{\Zp}{\GAL(K_\infty/K_0)}\hookrightarrow\power{\ValR}{\GAL(K_\infty/\Rat)}.
$$
The structure theorem for $\power{\Zp}{T}$-modules (\cite{Wa1},
Prop. 13.19) shows that $M_{(p)}=0$, since $\power{\Zp}{T}/(p^n)$
is not finitely generated over $\Zp$ for integers $n\geq0$. The
statement for $\power{\ValR}{\GAL(K_\infty/\Rat)}$ follows because
$i^{-1}(\ideal{p})=(p)$.
\end{proof}

We will now determine the group of invertible ideals of $\Omega$.
Since $\Omega$ is semilocal, it is given by
$$
\cart(\Omega)=\mult{\quotient{\Omega}}/\mult{\Omega}.
$$
In our main statement we compare two elements of $\cart(\Omega)$.
If $\ValR'$ is a faithfully flat extension of $\ValR$, e.g. the
valuation ring of a finite extension of $\quotient{\ValR}$, then
the induced map
$$
\cart(\Omega)\To\cart(\ValR'\tensorp{\ValR}\Omega)
$$
is injective.  Therefore, the above assumption that $\ValR$
contains the values of the characters of $\GAL(K_\infty/\Rat)$ is
no restriction for our purposes. (Alternatively, it can be
circumvented by using components instead of characters, as in
\cite{MW1}.)

If $\Rat_\infty\subset L_\infty \subset K_\infty$ is any
intermediate extension, we write
$$
\prKL[K_{\infty}/L_{\infty}]:\power{\ValR}{\GAL(K_{\infty}/\Rat)}\To\power{\ValR}{\GAL(L_\infty/\Rat)}
$$
for the natural projection. Note that the ring homomorphism
$\prKL[K_{\infty}/L_{\infty}]$ is extendable. Indeed, the induced
map $\widetilde{\prKL}$ between the normalisations of both rings
is extendable for almost trivial reasons. Since the inclusion
$\Omega\To\widetilde{\Omega}$ is extendable and maps zero divisors
to zero divisors, it follows that $\prKL$ is extendable as well.

\section{Iwasawa Cohomology}\label{sec:IwasawaCohomology}

Consider an open subscheme $U$ of $\SPEC \Int$ and let $S$ denote
its closed complement. If $F/\Rat$ is a finite field extension, we
set
$$
U_F=U\times\SPEC \Order[F],\qquad S_F=S\times\SPEC \Order[F],
$$
where $\Order[F]$ denotes the ring of integers of $F$. Write
$$
j_F: \SPEC F\To U_F
$$
for the inclusion of the generic point. As before, $K_\infty$ will
denote the cyclotomic $\Zp$-extension of an abelian number field.

Let $M(\rho)$ be a finitely generated $\ValR$-module $M$ together
with a continuous representation
$$
\rho: \GAL(\algclosure{\Rat}/\Rat)\To\AUT_{\ValR}M
$$
(where we give $Aut_{\ValR} M$ its profinite topology). Let
further
$$
\iota: \GAL(\algclosure{\Rat}/\Rat)\To\mult{\Omega}, \quad
\iota(g)=\bar{g}^{-1}\in\GAL(K_{\infty}/\Rat)
$$
denote the contragredient of the natural representation. The
$\algebra{\Omega}{\GAL(\algclosure{\Rat}/\Rat)}$-module
$$
\Ind M(\rho)=\Omega(\iota)\tensorp{\ValR}M(\rho)
$$
gives rise to a projective system of \'etale sheaves
$$
j_{\Rat*}\Ind M(\rho)=
\bigl(j_{\Rat*}\bigl(\algebra{\ValR/p^n\ValR}{\GAL(K_n/\Rat)}(\iota)\tensorp{\ValR}M(\rho)
\bigr)\bigr)_{n=1}^{\infty}
$$
on $U$. (We reemphasise that the action of
$\GAL(\algclosure{\Rat}/\Rat)$ on the module
$\algebra{\ValR/p^n\ValR}{\GAL(K_n/\Rat)}(\iota)$ is given by
$\iota$, i.e. $\algebra{\ValR/p^n\ValR}{\GAL(K_n/\Rat)}$ is
considered as trivial $\GAL(\algclosure{\Rat}/\Rat)$-module.)

\begin{defn}\label{defn:Iwasawa complex}
We define the Iwasawa complex of $M(\rho)$ over $U$ to be the
cochain complex of continuous \'etale cohomology
$$
\RightD\IwaSect(U_{K_{\infty}},M(\rho))=\RightD(\varprojlim_n\etSect)(U,j_{\Rat
*}\Ind M(\rho)),
$$
as constructed by U.~Jannsen in \cite{Ja1}. If $Z$ is a closed
subscheme of $U$, we define
$$
\RightD\IwaSect(U_{K_{\infty}},Z,M(\rho))=
\RightD(\varprojlim_n{\etSect})(U,Z,j_{\Rat *}\Ind M(\rho))
$$
to be the complex of continuous \'etale cohomology with support in
$Z$. These complexes are to be understood as objects of the
derived category of $\Omega$-modules. Their $i$-th cohomology
modules will be denoted by $\IwaH^i(U_{K_\infty},M(\rho))$,
respectively $\IwaH^i(U_{K_\infty},Z,M(\rho))$.
\end{defn}

\begin{rem}
Alternatively, it should also be possible to use the formalism of
T.~Ekedahl \cite{Ek}.
\end{rem}

Here are some basic properties of
$\RightD\IwaSect(U_{K_{\infty}},M(\rho))$.

\begin{prop}\label{prop:IwaCohomology}
Assume $p\notin U$.
\begin{enumerate}
\item For all $i\in\Int$,
$$
\IwaH^i(U_{K_\infty},M(\rho))= \varprojlim_n\etH^i(U_{K_n},j_{K_n
*}M(\rho)),
$$
where the limit is taken with respect to the corestriction maps.
\item In particular,
$$
\IwaH^0(U_{K_\infty},M(\rho))=0.
$$
\end{enumerate}
\end{prop}
\begin{proof}
By \cite{Mi2}, Theorem~II.2.13 the modules
$\etH^i(U_{K_n},j_{K_n*}M(\rho))$ are finite. The asserted
equality in (i) follows by \cite{Ja1}, Proposition~1.6 and
Lemma~1.15.

It remains to verify that $\IwaH^0(U_{K_\infty},M(\rho))=0$. As
$M$ is noetherian, there is an $n_0$ such that the inflation map
\begin{multline*}
\etH^0(U_{K_{n_0}},j_{K_{n_0}*}M(\rho))=\\
M(\rho)^{\GAL(\algclosure{\Rat}/K_{n_0})}\hookrightarrow M(\rho)^{\GAL(\algclosure{\Rat}/K_n)}\\
=\etH^0(U_{K_n},j_{K_n*}M(\rho))
\end{multline*}
is the identity for $n\geq n_0$ and therefore, the corestriction
map is multiplication by $p^{n-n_0}$. Hence, the limit over the
corestriction maps vanishes.
\end{proof}

If both $U_{K_{\infty}}$ and $\rho$ are unramified over $U$, then
all sheaves in the projective system $j_{\Rat*}\Ind M(\rho)$ are
locally constant. Under these circumstances one can identify
$\RightD\IwaSect(U_{K_{\infty}},M(\rho))$ with the complex of
continuous cochains of the topological $\fG(U)$-module $\Ind
M(\rho)$, where $\fG(U)$ denotes the \'etale fundamental group of
$U$ (see Proposition~II.2.9 of \cite{Mi2}). This setting has been
extensively explored by J.~Nekov{\'a}{\v{r}} in \cite{Nek1}. We
recall some of the consequences.

\begin{prop}\label{prop:cohom dim}
Assume that $U_{K_{\infty}}$ and $\rho$ are unramified over $U$.
Then $\RightD\IwaSect(U_{K_{\infty}},M(\rho))$ is acyclic outside
degrees $1$ and $2$.
\end{prop}
\begin{proof}
The cohomological $p$-dimension of $\fG(U)$ is $2$ (see
\cite{Ne1}, Theorem~8.3.19).
\end{proof}

\begin{prop}\label{prop:tensor products}
Let $M(\rho)$ be free as an $\ValR$-module and let $W$ be a
finitely generated $\Omega$-module. Assume that $U_{K_{\infty}}$
and $\rho$ are unramified over $U$. Then there exists a natural
quasi-isomorphism
$$
W\dtensorp{\Omega}\RightD\IwaSect(U_{K_\infty},M(\rho))=
\RightD(\varprojlim_n\etSect)(U,j_{\Rat *}(W\tensorp{\Omega}\Ind
M(\rho))).
$$
\end{prop}
\begin{proof}
See \cite{Nek1}, Proposition~3.4.4.
\end{proof}

Observe that for any
$\rho:\GAL(\algclosure{\Rat}/\Rat)\To\mult{\ValR}$ we can find an
abelian cyclotomic $\Zp$-extension $K_{\infty}$ such that $\rho$
factors through $\GAL(K_\infty/\Rat)$. We will then denote by
$$
\Tw{\rho}:\power{\ValR}{\GAL(K_{\infty}/\Rat)}\To\power{\ValR}{\GAL(K_\infty/\Rat)}
$$
the ring automorphism that maps $g\in\GAL(K_\infty/\Rat)$ to
$\rho(g)g$. For any ring homomorphism $f:R\To S$ and any
$R$-module $M$ we write
$$
f_*M=S\tensorp{R}M
$$
for the base extension to $S$.

\begin{prop}\label{prop:base change comps}
Let $M(\rho)$ be free as an $\ValR$-module. Assume that $\rho$ is
unramified outside a finite set of primes. Choose $U$ such that
for any $l\in U$ the ramification index of $l$ in
$K_{\infty}/\Rat$ is prime to $p$. Then
\begin{enumerate}
\item $\RightD\IwaSect(U_{K_\infty},M(\rho))$ is perfect,

\item for any intermediate field $\Rat_\infty\subset
L_\infty\subset K_\infty$ there exists a natural quasi-isomorphism
$$
\LeftD\prKL[K_{\infty}/L_{\infty}*]\RightD\IwaSect(U_{K_\infty},M(\rho))=
\RightD\IwaSect(U_{L_\infty},M(\rho)),
$$

\item For any $\chi:\GAL(K_\infty/\Rat)\To\mult{\ValR}$ there
exist a natural quasi-isomorphism
$$
\Tw{\chi*}\RightD\IwaSect(U_{K_\infty},M(\rho))=
\RightD\IwaSect(U_{K_\infty},M(\chi^{-1}\rho)).
$$
\end{enumerate}
\end{prop}
\begin{proof}
Let $V$ be an open subscheme of $U$ such that both $\rho$ and
$V_{K_\infty}$ are unramified over $U$. The localisation triangle
$$
\RightD\IwaSect(U_{K_\infty},U-V,M(\rho))\To
\RightD\IwaSect(U_{K_\infty},M(\rho))\To
\RightD\IwaSect(V_{K_\infty},M(\rho))
$$
(\cite{Ja1}, 3.6) implies that to prove (i), it is sufficient to
show that the two outer complexes are perfect. The right complex
is immediately seen to be perfect by Proposition~\ref{prop:tensor
products}. We will prove in Proposition~\ref{prop:local factors}
that the left complex is perfect as well.

By Remark \ref{rem:Local factors tensor} and the localisation
triangle it also suffices to prove (ii) and (iii) for the scheme
$V$. Claim (ii) then follows directly from the above proposition.
For (iii) it remains to notice that
\begin{align*}
\Tw{\chi*}\Ind M(\rho)&\To \Ind M(\chi^{-1}\rho), \\
1\tensor w\tensor m&\mapsto \Tw{\chi}(w)\tensor m\qquad (
w\in\Omega(\iota), m\in M(\rho))
\end{align*}
is a $\GAL(\algclosure{\Rat}/\Rat)$-equivariant isomorphism of
$\Omega$-modules.
\end{proof}

\section{Cohomology of
$\ValR(\cyclchar)$}\label{sec:Zp1Cohomology}

In this section we will calculate the cohomology of the
one-dimensional representation $\ValR(\cyclchar)$ given by the
cyclotomic character
$$
\cyclchar: \GAL(\algclosure{\Rat}/\Rat)\To\mult{\Zp}.
$$
The following proposition establishes the link to the objects of
classical Iwasawa theory.

\begin{prop}\label{prop:cohom of Zp(1)}
Let $U$ be an open subscheme of $X=\SPEC \Int$ such that $p$ lies
in the complement $S$ of $U$.
\begin{enumerate}
\item There exist a canonical isomorphism of $\Omega$-modules
$$
\IwaH^1(U_{K_\infty},\ValR(\cyclchar))=\varprojlim_n
\ValR\tensorp{\Int}\Gm(U_{K_n}),
$$
where $\Gm$ denotes the multiplicative group.
\item The following
sequence of $\Omega$-modules is exact
\begin{gather*}
0\To\varprojlim_n\ValR\tensorp{\Int}\Gm(X_{K_n})\To
\IwaH^1(U_{K_\infty},\ValR(\cyclchar))\\
\To\varprojlim_n\etH^0(S_{K_n},\ValR)
\To\varprojlim_n\ValR\tensorp{\Int}\Pic(X_{K_n})\To\\
\IwaH^2(U_{K_\infty},\ValR(\cyclchar))
\To\varprojlim_n\etH^1(S_{K_n},\ValR)\To\ValR\To 0.
\end{gather*}
\end{enumerate}
\end{prop}
\begin{proof}
This is proved in the same way as \cite{BG1}, Proposition~5.1. The
idea is to combine the calculation of the cohomology groups of
$\Gm$ in \cite{Mi2}, Proposition~II.2.1, with the Kummer exact
sequence on $U_{K_n}$ and then to pass to the limit. The last term
can then be identified as the tensor product of
$$
\Zp=\varprojlim_n \KER
(\etH^3(X_{K_n},\Gm)\xrightarrow{p^n}\etH^3(X_{K_n},\Gm))
$$
with $\ValR$.
\end{proof}

This result is complemented by the following

\begin{prop}\label{prop:gysin factors}
There exist (non-canonical) isomorphisms of $\Omega$-modules
\begin{align*}
\varprojlim_n\etH^0(S_{K_n},\ValR)
&\isomorph\algebra{\ValR}{\GAL(K_{\infty}/\Rat)/D_p}\\
\varprojlim_n\etH^1(S_{K_n},\ValR) &\isomorph\Dsum_{l\in
S}\algebra{\ValR}{\GAL(K_{\infty}/\Rat)/D_l},
\end{align*}
where $D_l$ denotes the decomposition subgroup of the prime $l$ in
$\GAL(K_{\infty}/\Rat)$.
\end{prop}
\begin{proof}
One of the fundamental properties of $\Zp$-extensions is the fact
that $K_\infty/K_0$ is unramified outside the primes over $p$ (see
\cite{Wa1}, Prop. 13.2). For cyclotomic $\Zp$-extensions one also
knows that there exists a number $n_0$ such that all primes over
$p$ are totally ramified in $K_{\infty}/K_{n_0}$ and such that
none of the primes in $S_{K_{n_0}}$ splits in $K_{\infty}/K_{n_0}$
(see \cite{Wa1}, Ex. 13.2). In particular, $S_{K_\infty}\To
S_{K_n}$ is a homeomorphism for $n\geq n_0$.

On the other hand,
$$
\etH^0(S_{K_n},\ValR)\isomorph\Dsum_{v\in
S_{K_n}}\ValR\isomorph\etH^1(S_{K_n},\ValR).
$$
An elementary calculation shows that the corestriction map
$$
\etH^i(S_{K_{n+1}},\ValR)\To\etH^i(S_{K_n},\ValR)
$$
for $n\geq n_0$ is the identity for $i=1$, respectively the
multiplication by the residue degree of $v$ on the $v$-component
of $\etH^0(S_{K_{n+1}},\ValR)$ for $i=0$. But the residue degree
is $1$ or $p$ depending on wether $v$ lies over $p$ or not. Now
pass to the limit. The choice of an element of $S_{K_\infty}$ for
each prime in $S$ induces the desired isomorphisms.
\end{proof}

\begin{cor}\label{cor:local factors Zp(1)}
Let $T$ be a closed subscheme of $U$. Then the complex
$\RightD\IwaSect(U_{K_\infty},T,\ValR(\cyclchar))$ is acyclic
outside degree $3$ and
$$
\IwaH^3(U_{K_{\infty}},T,\ValR(\cyclchar))\isomorph\Dsum_{l\in
T}\algebra{\ValR}{\GAL(K_{\infty}/\Rat)/D_l}
$$
\end{cor}
\begin{proof}
Easy application of the snake lemma.
\end{proof}

\section{Local Factors}\label{sec:LocalFactors}

In this section we shall examine the relative cohomology complexes
$\RightD\IwaSect(U_{K_\infty},S,M(\rho))$ for arbitrary continuous
$\GAL(\algclosure{\Rat}/\Rat)$-representations $M(\rho)$, where
$M$ is any finitely generated $\ValR$-module. This will also
complete the proof of Proposition~\ref{prop:base change comps}.
Our aim is to extend Corollary \ref{cor:local factors Zp(1)} as
follows.

\begin{prop}\label{prop:local factors}
Let $U$ be an open subscheme of $\SPEC \algebra{\Int}{1/p}$, $S$ a
closed subscheme of $U$. Then
\begin{enumerate}
\item $\RightD\IwaSect(U_{K_\infty},S,M(\rho))$ is acyclic outside
degree $3$;

\item $\IwaH^3(U_{K_\infty},S,M(\rho))$ is a finitely generated
$\ValR$-module;

\item if for all $l\in S$ the prime $p$ does not divide the
ramification index of $l$ in $K_\infty/\Rat$, then
$\RightD\IwaSect(U_{K_\infty},S,M(\rho))$ is a perfect torsion
complex of $\power{\ValR}{\GAL(K_{\infty}/\Rat)}$-modules.
\end{enumerate}
\emph{Supplement:} If $K_{\infty}/\Rat$ is a $p$-extension,
$M=\ValR$, and $\rho$ is ramified over a prime $l$ that is
unramified in $K_{\infty}/\Rat$, then
$$
\CHAR \RightD\IwaSect(U_{K_\infty},l,\ValR(\rho))=(1).
$$
\end{prop}

\begin{proof}
As in the proof of Proposition~\ref{prop:gysin factors}, we fix a
number $m$ such that none of the primes in $S_{K_{m}}$ splits or
ramifies in $K_{\infty}/K_{m}$.

Write $v_m\in S_{K_m}$ for the image of a point $v\in
S_{K_\infty}$. From \cite{Ja1}, Prop. 3.8, we obtain an
isomorphism of $\power{\ValR}{\GAL(K_\infty/K_{m})}$-modules
\begin{multline*}
\IwaH^i(U_{K_\infty},S,M(\rho))=\\
\Dsum_{v\in S_{K_\infty}}\RightD^i(\varprojlim\etSect)(\SPEC
\Order[v_m]^h,v_m,j_{K_{m} *}\Ind[K_{\infty}/K_{m}]M(\rho)),
\end{multline*}
where $\Order[v_m]^h$ is the henselisation of the local ring at
$v_m$ and the limit is taken over the projective system $j_{K_{m}
*}\Ind[K_{\infty}/K_{m}]M(\rho)$.

We will now use the connection between \'{e}tale and Galois
cohomology. Fix a $v\in S_{K_\infty}$ and let $G_{v_m}$ denote the
absolute Galois group of $\quotient{\Order[v_m]^h}$, $I_{v_m}$ its
inertia subgroup and $g_{v_m}=G_{v_m}/I_{v_m}$ the Galois group of
the residue field of $v_m$. For any profinite group $G$, let
$\cd_p G$ denote the cohomological $p$-dimension, i.e. the largest
number $i$ such that the $i$-th group cohomology functor is
non-trivial on finite $p$-torsion $G$-modules. We have
\begin{align*}
\cd_p G_{v_m}&=2,\\
\cd_p I_{v_m}=\cd_p g_{v_m}&=1
\end{align*}
(see \cite{Ne1}, Prop. 3.3.4, Prop. 7.1.8). Further, it is well
known that for these groups, the cohomology groups of finite
$p$-torsion modules will again be finite. In particular, we may
interchange projective limits and continuous cohomology functors
during the subsequent considerations.

By \cite{Mi2}, Proposition~II.1.1.(b) and \cite{Mi1}, Ex. II.3.15
it follows that
\begin{multline*}
\RightD^i(\varprojlim\etSect)(\SPEC \Order[v_m]^h,j_{K_{m}
*}\Ind[K_{\infty}/K_{m}]M(\rho))=\\
\Hc^i(g_{v_m},H^0(I_{v_m},\Ind[K_{\infty}/K_{m}]M(\rho))).
\end{multline*}
Comparing the localisation sequence for
$$
\SPEC \quotient{\Order[v_m]^h}\hookrightarrow\SPEC
\Order[v_m]^h\hookleftarrow v_m
$$
with the Hochschild-Serre spectral sequence for $I_{v_m}\subset
G_{v_m}$ we obtain from the above
\begin{multline*}
\RightD^i(\varprojlim\etSect)(\SPEC \Order[v_m]^h,v_m,j_{K_{m}
*}\Ind[K_{\infty}/K_{m}]M(\rho))=\\
\Hc^{i-2}(g_{v_m},\Hc^1(I_{v_m},\Ind[K_{\infty}/K_{m}]M(\rho))).
\end{multline*}
As $v_m$ is unramified in $K_{\infty}/K_{m}$, the action of
$I_{v_m}$ on $\power{\ValR}{\GAL(K_{\infty}/K_{m})}(\iota)$ is
trivial and hence,
$$
\Hc^1(I_{v_m},\Ind[K_{\infty}/K_{m}]M(\rho))=
\Ind[K_{\infty}/K_{m}]\Hc^1(I_{v_m},M(\rho)).
$$
Observe that $\Hc^1(I_{v_m},M(\rho))$ is a finitely generated
$\ValR$-module.

To finish the proof of Proposition~\ref{prop:local factors}, (i)
and (ii), it remains to verify the following

\begin{lem}
Let $N(\tau)$ be a finitely generated $\ValR$-module $N$ together
with a continuous representation $\tau:g_{v_m}\To\AUT_{\ValR}N$.
Then
\begin{enumerate}
\item $\Hc^0(g_{v_m},\Ind[K_\infty/K_{m}]N(\tau))=0$,

\item $\Hc^1(g_{v_m},\Ind[K_\infty/K_{m}]N(\tau))$ is a finitely
generated $\ValR$-module.
\end{enumerate}
\end{lem}
\begin{proof}
By our assumption on $m$ we have $I_{v_n}=I_{v_m}$ and
$G_{v_n}/G_{v_m}=\GAL(K_n/K_{m})$ for $n\geq m$. Thus,
$$
\Hc^i(g_{v_m},\Ind[K_{\infty}/K_{m}]N(\tau))=\varprojlim_n\Hc^i(g_{v_n},N(\tau)).
$$
The same argument as in Proposition~\ref{prop:IwaCohomology}.(ii)
implies that this term vanishes for $i=0$. This proves (i). Claim
(ii) follows because $\Hc^1(g_{v_m},N(\tau))$ is a quotient of $N$
for all $n$.
\end{proof}

We now prove Proposition~\ref{prop:local factors}, (iii). After
decomposing by characters we may assume that $K_{\infty}$ is a
$p$-extension. In particular, $l\in S$ is unramified. By the same
argument as above, replacing $K_{m}$ by $\Rat$, we obtain
$$
\IwaH^3(U_{K_\infty},l,M(\rho))= \Hc^1(g_l,\Hc^1(I_l,\Ind
M(\rho))).
$$
Write $\Omega=\power{\ValR}{\GAL(K_{\infty}/\Rat)}$ and note that
$I_l$ acts trivially on $\Omega(\iota)$. Consequently,
$$
\Hc^1(I_l,\Ind M(\rho))=\Ind N(\tau),
$$
where $N=\Hc^1(I_l,M(\rho))$ is a finitely generated
$\ValR$-module and $\tau:g_l\to\AUT_{\ValR}N$ is the induced
representation. Recall that $g_l$ is topologically generated by
the geometric Frobenius element $\geoF_l$. By (i) the sequence
$$
0\To \Omega\tensorp{\ValR}N  \xrightarrow{\id -
\iota(\geoF_l)\tensor\tau(\geoF_l)}\Omega\tensorp{\ValR}N\To\IwaH^3(U_{K_\infty},l,M(\rho))\To
0
$$
is exact. As $\ValR$ is regular, $\Omega\tensorp{\ValR}N$ is
perfect as $\Omega$-complex; hence, so is
$\IwaH^3(U_{K_\infty},l,M(\rho))$ . The latter module is
$\Omega$-torsion, because it is a finitely generated
$\ValR$-module. This proves (iii).

To prove the supplement it suffices to recall that
$$
N=(M(\rho)^{R_l})_{T_l},
$$
where $R_l$ is the ramification subgroup and $T_l=I_l/R_l$. If
$M=\ValR$ and the restriction of $\rho$ to $I_l$ is non-trivial,
then this module is clearly $\ValR$-torsion. By
Proposition~\ref{prop:properties of char},(iii) and (v), applied
to the above exact sequence, we obtain
$$
\CHAR \RightD\IwaSect(U_{K_\infty},l,M(\rho))=\CHAR
(\Omega\tensorp{\ValR}N)\CHAR^{-1}(\Omega\tensorp{\ValR}N)=(1).
$$
This finishes the proof of Proposition~\ref{prop:local factors}.
\end{proof}

\begin{rem}\label{rem:Local factors tensor}
Let $W$ be a finitely generated $\Omega$-module. If either $W$ or
$\Hc^1(I_l,M(\rho))$ is flat as an $\ValR$-module, then
\begin{multline*}
W\dtensorp{\Omega}\RightD\IwaSect(U_{K_\infty},l,\Ind[K_{\infty}/\Rat]M(\rho))\isomorph\\
\Hc^1(g_l,W\tensorp{\Omega}\Ind[K_{\infty}/\Rat]\Hc^1(I_l,M(\rho)))\\
\isomorph\RightD\IwaSect(U_{K_\infty},l,W\tensorp{\Omega}\Ind[K_{\infty}/\Rat]M(\rho)).
\end{multline*}
This is not true without the additional flatness assumption.
\end{rem}

\section{Cyclotomic Elements and
$L$-Elements}\label{sec:CyclotomicElements}

The aim this section is to assign to each admissible triple
$(K_\infty,\rho,U)$ given by
\begin{itemize}
    \item a cyclotomic $\Zp$-extension $K_\infty$ of an abelian number
    field,
    \item a representation
    $\rho:\GAL(\algclosure{\Rat}/\Rat)\To\mult{\ValR}$,
    \item an open subscheme $U$ of $\SPEC \Int$ that does not
    contain any prime whose ramification index in $K_\infty/\Rat$
    is divisible by $p$
\end{itemize}
an $L$-element
$L(U_{K_\infty},\rho^{-1}\cyclchar)\in\mult{\quotient{\power{\ValR}{\GAL(K_\infty/\Rat)}}}$
and a cyclotomic element $c(U_{K_\infty},\rho)\in
\IwaH^1(U_{K_\infty},\ValR(\rho))$.

If $\rho$ and $U_{K_\infty}$ are unramified over $U$, then our
definition follows along the lines of the classical construction
(see \cite{Wa1}, \S7.2, respectively \cite{Ru1}, \S3.2). By the
Theorem of Kronecker and Weber there exists a number $f$ such that
\begin{itemize}
    \item the set of prime divisors of $f$ is equal to the
    complement of $U$ in $\SPEC \algebra{\Int}{1/p}$,
    \item $K_\infty\subset \CyclF{fp^{\infty}}$,
    \item $\rho$ factors through $\GAL(\CyclF{fp^{\infty}}/\Rat)$.
\end{itemize}
Let $\geoF_l$ denote the geometric Frobenius element and set
$$
\geoF_a=\prod_{l \text{ prime}}\geoF_l^{v_l(a)}
$$
for each positive integer $a$. The Stickelberger elements
$$
\Stickel{fp^k}=\sum_{\substack{0<a<fp^k\\
\ggt{a}{fp}=1}}(\frac{a}{fp^k}-\frac{1}{2})\geoF_{a}\in\algebra{\Rat_p}{\GAL(\CyclF{fp^k}/\Rat)}
$$
are compatible under the projection maps induced by
$$
\GAL(\CyclF{fp^{k+1}}/\Rat)\To\GAL(\CyclF{fp^k}/\Rat)
$$
and define an element
$$
\Stickel{fp^{\infty}}\in\quotient{\power{\ValR}{\GAL(\CyclF{fp^{\infty}}/\Rat)}}.
$$
Further, we fix for each number $k$ a primitive root of unity
$\pru{k}$ such that $\pru{ks}^s=\pru{k}$. By
Proposition~\ref{prop:cohom of Zp(1)} we may regard the system
$c_{fp^{\infty}}=(1-\pru{fp^k})_{k=0}^{\infty}$ as an element of
$\IwaH^1(U_{\CyclF{fp^{\infty}}},\ValR(\cyclchar))$.

\begin{defn}
Let $\rho$ and $U_{K_\infty}$ be unramified over $U$. Denote by
$\pr{+}$ and $\pr{-}$ the projectors onto the $(+1)$-eigenspace,
respectively the $(-1)$-eigenspace of the complex conjugation and
chose $f$ as above. We set
$$
L(U_{K_\infty},\rho^{-1}\cyclchar)=
\prKL[\CyclF{fp^{\infty}}/K_\infty]\Tw{\rho^{-1}\cyclchar}(\pr{+}-\pr{-}\Stickel{fp^\infty})
$$
The cyclotomic element $c(U_{K_\infty},\rho)$ is defined to be the
image of $1\tensor\pr{+}c_{fp^{\infty}}$ under the homomorphism
$$
(\prKL[\CyclF{fp^{\infty}}/K_\infty]\Tw{\rho^{-1}\cyclchar})_*
\IwaH^1(U_{\CyclF{fp^{\infty}}},\ValR(\cyclchar))\To
\IwaH^1(U_{K_{\infty}},\ValR(\rho)).
$$
For any $l\in U$ we denote by
$$
E_l(K_\infty,\rho^{-1}\cyclchar)=
1-\prKL[\CyclF{fp^{\infty}}/K_\infty]\Tw{\rho^{-1}\cyclchar}(\geoF_l)
$$
the Euler factor at $l$.
\end{defn}

We will now extend this definition to arbitrary admissible triples
$(K_{\infty},\rho,U)$. Let $P_{\infty}$ be the maximal
$p$-extension inside $K_\infty$. We may decompose
$\power{\ValR}{\GAL(K_\infty/\Rat)}$ by the characters of
$\GAL(K_\infty/P_\infty)$. The $L$-elements and cyclotomic
elements for $K_\infty$ are then completely determined by their
projections onto the components of the corresponding decomposition
of $\mult{\quotient{\power{\ValR}{\GAL(K_\infty/\Rat)}}}$,
respectively $\IwaH^1(U_{K_\infty},\ValR(\rho))$. Consequently, it
suffices to consider triples $(P_{\infty},\rho,U)$ with
$P_{\infty}$ a $p$-extension. Note that every ramification index
is now a power of $p$. By assumption, $U_{P_\infty}$ is therefore
an unramified cover of $U$. In other words, we only need to deal
with the ramification of $\rho$.

Let $V\subset U$ be the maximal open subscheme such that $\rho$ is
unramified over $V$ (note that the complement is a finite set
because $\rho$ is one-dimensional ) and set
\begin{align*}
L(U_{P_\infty},\rho^{-1}\cyclchar)&=L(V_{P_\infty},\rho^{-1}\cyclchar),\\
c(U_{P_\infty},\rho)&=c(V_{P_\infty},\rho),\\
\intertext{and for any $l\in U-V$}
E_l(P_\infty,\rho^{-1}\cyclchar)=1.
\end{align*}

It is easy to check that this definition is consistent with the
previous construction. Moreover, we have

\begin{prop}\label{prop:compatibilities}
Let $(K_{\infty},U,\rho)$ be any admissible triple.
\begin{enumerate}
\item The elements $L(U_{K_\infty},\rho^{-1}\cyclchar)$,
$c(U_{K_\infty},\rho)$, and $E_l(K_{\infty},\rho^{-1}\cyclchar)$
are compatible under the projection maps
$\prKL[K_{\infty}/L_\infty]$ and under twists by continuous
characters $\GAL(K_{\infty}/\Rat)\To\mult{\ValR}$.

\item Let $V\subset U$ be an open subscheme with closed complement
$T=U-V$. Then
\begin{align*}
L(V_{K_{\infty}},\rho)&=
\Tw{\rho}(\pr{+})+\Tw{\rho}(\pr{-})L(U_{K_{\infty}},\rho)
\prod_{l\in T}E_l(K_{\infty},\rho),\\
c(V_{K_\infty},\rho)&= c(U_{K_\infty},\rho)\prod_{l\in
T}E_l(K_{\infty},\rho^{-1}\cyclchar).
\end{align*}
\end{enumerate}
\end{prop}
\begin{proof}
This is implied by the corresponding compatibility properties of
$\Stickel{fp^{\infty}}$ and $c_{fp^{\infty}}$, respectively true
by definition.
\end{proof}

We want to sketch briefly the connection between our $L$-elements
and the Kubota-Leopoldt $L$-function $L_p(s,\chi)$. Let
$K_\infty=\Rat_\infty$ and $U=\SPEC \algebra{\Int}{1/p}$.
According to the decomposition $\mult{\ValR}=\mu\times F$ into the
torsion group of roots of unity $\mu$ and the torsion-free
$\Zp$-module $F$ we can write
$$
\rho=\rho_f\rho_\infty,
$$
where $\rho_f$ takes values in $\mu$ and
$$
\rho_\infty:\GAL(\Rat_\infty/\Rat)\To F
$$
is a continuous group homomorphism.

Let $N$ be the conductor of $\rho_f$ and
$$
\rec:\GAL(\CyclF{N}/\Rat)\To\mult{(\Int/N\Int)}
$$
the isomorphism which maps the geometric Frobenius $\geoF_l$ to
$l$, when $l$ is prime to $N$. Further, note that the number
$$
s=\frac{\log_p\rho_\infty(\gamma)}{\log_p\big((\cyclchar)_{\infty}(\gamma)\big)}\in
\log_p F\subset\Cp
$$
does not depend on the choice of a nontrivial
$\gamma\in\GAL(\Rat_\infty/\Rat)$. Hence, we may write
$\rho_\infty=(\cyclchar)_\infty^s$.

\begin{prop}\label{prop:evaluation of L-element}
Let $\rho=(\cyclchar)_{\infty}^s\rho_f$ be an even one-dimensional
representation. Then
$$
\prKL[\Rat_\infty/\Rat](L(U_{\Rat_\infty},\cyclchar\rho))=L_p(1+s,\rho_f\conjunct\rec^{-1}).
$$
\end{prop}
\begin{proof}
See \cite{Wa1}, Theorem~7.10, but observe that the identification
$$
\GAL(\CyclF{N}/\Rat)\To\mult{(\Int/N\Int)}
$$
used in loc. cit. is given by $1/\rec$. In particular,
$w\conjunct\rec=(\cyclchar)_f^{-1}$, where $w$ denotes the
Teichm\"uller character.
\end{proof}

\begin{rem}\label{rem:relation to f}
Let $N$ be prime to $p$ and let
$\chi:\mult{(\Int/Np\Int)}\To\mult{\ValR}$ be an odd Dirichlet
character of conductor $N$ or $Np$ (i.e. $\chi$ is of the first
kind). Set $U=\SPEC\algebra{\Int}{1/p}$ and let
$$
\beta:\power{\ValR}{\GAL(\Rat_\infty/\Rat)}\To\power{\ValR}{T}
$$
be the isomorphism that maps $\geoF_{1+Np}^{-1}$ to $T+1$. By
construction we then have
$$
f(T,\chi\omega)=\beta(L(U_{\Rat_{\infty}},\chi\conjunct\rec)),
$$
where $f(T,\chi\omega)$ is the element introduced in \cite{Wa1},
\S7.2.
\end{rem}

\begin{rem}
Let $\chi$ be a finite character and $k$ an integer.
Proposition~\ref{prop:evaluation of L-element} implies that our
element $L(U_{\Rat_\infty},\cyclchar^{1-k}\chi)$, with $U=\SPEC
\algebra{\Int}{1/p}$, coincides with the $p$-adic $L$-function
$\mathcal{L}_p(\chi,1-k)$ used in \cite{HK1}. However, note that
there is a sign error in the definition of this function. The
correct definition should read as follows (in the notation of loc.
cit.). For all $\widetilde{\mathcal{O}}_p$ and all characters
$\tau: \Gamma\To\widetilde{\mathcal{O}}_p^*$ of finite order,
$$
\tau(\mathcal{L}_p(\chi,1-k))=(1-\chi\tau(p)p^{k-1})L(\chi\tau,1-k).
$$
\end{rem}

\begin{rem}
The elements $c(U_{K_\infty},\rho)$ depend on the choice of the
system of roots of unity $(\pru{k})$, but the submodule of
$\IwaH^1(U_{K_\infty},\ValR(\rho))$ generated by
$c(U_{K_\infty},\rho)$ does not. This is the actual object we are
interested in.
\end{rem}

The Theorem of Ferrero-Washington can be rephrased to the
statement that the $\mu$-invariants of the $p$-adic $L$-functions
vanish. The following proposition translates this formulation to
our setting.

\begin{prop}\label{prop:mu=0 for L-elements}
Let $(K_{\infty},U,\rho)$ be an admissible triple and
$l\in\SPEC\Int$ of prime-to-$p$ ramification in $K_{\infty}/\Rat$.
Then $L(U_{K_\infty},\rho)$ and $E_l(K_\infty,\rho)$ map to units
in $\power{\ValR}{\GAL(K_\infty/\Rat)}_{\ideal{p}}$ for each prime
ideal $\ideal{p}$ of codimension $1$ with $p\in\ideal{p}$.
\end{prop}
\begin{proof}
By Proposition~\ref{prop:compatibilities}.(iii) we may assume that
$U_{K_\infty}$ and $\rho$ are unramified over $U$. Choose $f$ as
above and observe that
$$
(1-(1+fp)\geoF_{1+fp})\Stickel{fp^{\infty}}
\in\power{\ValR}{\GAL(\CyclF{fp^{\infty}}/\Rat)}.
$$
Define
$$
h(U_{K_\infty},\rho)=1-(1+fp)\prKL[\CyclF{fp^{\infty}}/K_\infty]\Tw{\rho}(\geoF_{1+fp}).
$$
We need to show that neither of $h(U_{K_\infty},\rho)$,
$hL(U_{K_\infty},\rho)$, or $E_l(K_\infty,\rho)$ is contained in
$\ideal{p}$. For this, we can replace
$\power{\ValR}{\GAL(K_\infty/\Rat)}$ by its normalisation and then
decompose by the characters of $\GAL(K_0/\Rat)$. Hence, we may
assume that $K_\infty=\Rat_\infty$. In particular, $\ideal{p}$ is
the radical of $(p)$. After twisting by an appropriate character
of $\GAL(\Rat_\infty/\Rat)$ we may further require that
$\rho=\chi\conjunct\rec^{-1}$ for a Dirichlet character $\chi$ of
the first kind.

Obviously, $E_l(\Rat_\infty,\chi)$ and $h(U_{\Rat_\infty},\chi)$
are prime to $p$ (note that the images of $\geoF_l$ and
$\geoF_{1+fp}$ are nontrivial in $\GAL(\Rat_\infty/\Rat)$). If
$\chi$ is even, then $L(U_{K_\infty},\chi)=1$. If $\chi$ is odd,
the claim for $hL(U_{K_\infty},\theta)$ is by Remark
\ref{rem:relation to f} equivalent to the vanishing of the
$\mu$-invariant of $f(T,\chi\omega)$, hence to the Theorem of
Ferrero-Washington (\cite{Wa1}, \S7.5, respectively \S16.2).
\end{proof}

\section{The Main Theorem} \label{sec:main}

Let $(K_\infty,\rho,U)$ be an admissible triple in the sense of
Section~\ref{sec:CyclotomicElements} and set
$$
\Omega=\power{\ValR}{\GAL(K_\infty/\Rat)},
$$
where $\ValR$ is the valuation ring of a finite extension of
$\Rat_p$. As explained in Section~\ref{sec:ProfGroupRing}, we may
assume without loss of generality that $\ValR$ contains the values
of all characters of $\GAL(K_0/\Rat)$.

Before we state our main theorem we will explain how to modify the
Iwasawa complex $\RightD\IwaSect(U_{K_\infty},\ValR(\rho))$ (see
Definition \ref{defn:Iwasawa complex}) by the cyclotomic element
$$
c(U_{K_\infty},\rho)\in\IwaH^1(U_{K_\infty},\ValR(\rho))
$$
introduced in the preceding section.

Recall that $\RightD\IwaSect(U_{K_\infty},\ValR(\rho))$ is acyclic
in degree $0$ (see Proposition~\ref{prop:IwaCohomology}.(ii)). In
particular, there exists a unique morphism
$$
\Omega
c(U_{K_\infty},\rho)[-1]\To\RightD\IwaSect(U_{K_\infty},\ValR(\rho))
$$
in the derived category that induces the natural inclusion on
cohomology.

\begin{defn}
Denote by $\Rc(U_{K_\infty},\ValR(\rho))$ the complex (unique up
to quasi-isomorphism) fitting into the following distinguished
triangle
$$
\Omega
c(U_{K_\infty},\rho)[-1]\To\RightD\IwaSect(U_{K_\infty},\ValR(\rho))\To\Rc(U_{K_\infty},\ValR(\rho)).
$$
\end{defn}

\begin{lem}\label{lem:torsion+free}\
\begin{enumerate}
\item $\Rc(U_{K_\infty},\ValR(\rho))$ is a perfect torsion complex
of $\Omega$-modules.

\item If $K_\infty$ is totally real and $\rho$ is odd, then
$\Omega c(U_{K_\infty},\rho)$ is a free $\Omega$-module of rank
$1$.
\end{enumerate}
\end{lem}
\begin{proof}
In Proposition~\ref{prop:base change comps}.(i) we have already
confirmed that the complex
$\RightD\IwaSect(U_{K_\infty},\ValR(\rho))$ is perfect. Hence, it
suffices to prove (ii) and that the complex in (i) is torsion.

Let $\widetilde{\Omega}$ denote the normalisation of $\Omega$ in
its total quotient ring $\quotient{\Omega}$. On the one hand we
have
$$
\quotient{\Omega}\tensorp{\Omega}\IwaH^i(U_{K_\infty},\ValR(\rho))=
\quotient{\Omega}\tensorp{\widetilde{\Omega}}
\Big(\Dsum_{\chi}\IwaH^i(U_{\Rat_\infty},\ValR(\chi^{-1}\rho))\Big),
$$
where the sum runs over all characters $\chi$ of $\GAL(K_0/\Rat)$;
on the other hand $\Omega
c(U_{K_\infty},\rho)\subset\IwaH^1(U_{K_\infty},\ValR(\rho))$ is a
free $\Omega$-module of rank $1$ if and only if
$\quotient{\Omega}c(U_{K_\infty},\rho)$ is a free
$\quotient{\Omega}$-module of rank $1$. This module decomposes as
$$
\quotient{\Omega}c(U_{K_\infty},\rho)=
\Dsum_{\chi}\quotient{\power{\ValR}{\GAL(\Rat_\infty/\Rat)}}c(U_{\Rat_\infty},\chi^{-1}\rho).
$$
Observe hereby that $\rho$ and $\chi^{-1}\rho$ have the same
parity if $K_\infty$ is totally real. Hence, it is enough to
consider the case $K_\infty=\Rat_\infty$. By
Proposition~\ref{prop:local factors} and
Proposition~\ref{prop:compatibilities} we may replace $U$ by
$\SPEC \algebra{\Int}{1/p}$, noting that the Euler factors
$E_l(\Rat_{\infty},\rho)$ are nonzero divisors. We are now reduced
to the statement of \cite{HK1}, Proposition~4.2.1. Observe that
the proof of (ii) in this situation uses K.~Kato's explicit
reciprocity law as an essential ingredient.
\end{proof}

\begin{cor}\label{cor:modified_compatibility}
The compatibility properties of Proposition~\ref{prop:base change
comps} hold for $\ValR(\rho)$, with $\RightD\IwaSect$ replaced by
$\Rc$.
\end{cor}
\begin{proof}
Easy consequence of Proposition~\ref{prop:compatibilities} and the
above lemma.
\end{proof}

We are now ready to formulate and prove our main result.

\begin{thm}\label{thm:Main Theorem}
Let $p$ be an odd prime, $K_\infty$ the cyclotomic $\Zp$-extension
of an abelian number field, and $U$ an open subscheme of $\SPEC
\Int$ such that the ramification index in $K_{\infty}/\Rat$ of
every place in $U$ is prime to $p$. Then
\begin{enumerate}
\item $\Rc(U_{K_\infty},\ValR(\rho))_{\ideal{p}}$ is acyclic for
all primes $\ideal{p}$ of codimension~$1$ of
$\power{\ValR}{\GAL(K_\infty/\Rat)}$ that contain $p$.

\item The characteristic ideal of $\Rc(U_{K_\infty},\ValR(\rho))$
is generated by the $L$-element
$L(U_{K_\infty},\rho^{-1}\cyclchar)$.
\end{enumerate}
\end{thm}
\begin{proof}
By the subsequent lemma we are allowed to enlarge or shrink the
scheme $U$ at our discretion.

\begin{lem}\label{lem:indep of Euler-f}
Let $V\subset U$ be an open subscheme of $U$. Then both statements
of Theorem~\ref{thm:Main Theorem} hold for $U$ if and only if they
holds for $V$.
\end{lem}

\begin{proof}
Set $T=U-V$ and define $\cmplx{C}$ by the following distinguished
triangle
$$
\Omega
c(V_{K_\infty},\rho)[-1]\To\RightD\IwaSect(U_{K_\infty},\ValR(\rho))\To\cmplx{C}.
$$
where the first map is induced by the inclusion
$$
\Omega
c(V_{K_\infty},\rho)\hookrightarrow\IwaH^1(V_{K_\infty},\ValR(\rho))
=\IwaH^1(U_{K_\infty},\ValR(\rho)).
$$

We obtain the following two triangles of perfect torsion complexes
\begin{gather*}
\RightD\IwaSect(U_{K_\infty},T,\ValR(\rho))\To\cmplx{C}\To\Rc(V_{K_\infty},\ValR(\rho))\\
\cmplx{C}\To\Rc(U_{K_\infty},\ValR(\rho))\To \Omega
c(U_{K_\infty},\rho)/\Omega c(V_{K_\infty},\rho).
\end{gather*}
Proposition~\ref{prop:compatibilities} implies
$$
\Omega c(U_{K_\infty},\rho)/\Omega c(V_{K_\infty},\rho)=
\Tw{\rho^{-1}\cyclchar}(\pr{+})\Big(\Omega\Big/\prod_{l\in
T}E_l(K_\infty,\rho^{-1}\cyclchar)\Omega\Big)
$$
and since $E_l(K_\infty,\rho^{-1}\cyclchar)$ is a unit of
$\Omega_\ideal{p}$ for any prime ideal $\ideal{p}$ of codimension
$1$ with $p\in\ideal{p}$ (see Proposition~\ref{prop:mu=0 for
L-elements}) it follows that
$$
\big(\Omega c(U_{K_\infty},\rho)/\Omega
c(V_{K_\infty},\rho)\big)_{\ideal{p}}=0.
$$
On the other hand, we know that
$\RightD\IwaSect(U_{K_\infty},T,\ValR(\rho))$ is acyclic outside
degree $3$ and that $\IwaH^3(U_{K_\infty},T,\ValR(\rho))$ is
finitely generated as $\ValR$-module (see
Proposition~\ref{prop:local factors}). Therefore,
$$
\RightD\IwaSect(U_{K_\infty},T,\ValR(\rho))_{\ideal{p}}\isomorph 0
$$
by Lemma~\ref{lem:criterion for vanishing}. This implies the
equivalence for part (i) of Theorem~\ref{thm:Main Theorem}.

By using the multiplicativity of the characteristic ideal, the
equivalence for part (ii) is reduced to proving that
$$
\CHAR
\IwaH^3(U_{K_\infty},l,\ValR(\rho))=E_l(K_\infty,\rho^{-1}\cyclchar)\Omega
$$
for $l\in T$. After decomposing by characters we may assume that
$K_\infty$ is a $p$-extension. For those primes over which $\rho$
is ramified the equality is implied by the supplement to
Proposition~\ref{prop:local factors}. For the remaining primes
choose $f$ as in Section~\ref{sec:CyclotomicElements}. From
Corollary \ref{cor:local factors Zp(1)} we obtain
\begin{multline*}
\CHAR \IwaH^3(U_{K_\infty},l,\ValR(\rho))=\\
\CHAR
(\prKL[\CyclF{fp^{\infty}}/K_{\infty}]\Tw{\rho^{-1}\cyclchar})_*
\algebra{\ValR}{\GAL(\CyclF{fp^{\infty}}/\Rat)/D_l}\\
=E_l(K_\infty,\rho^{-1}\cyclchar)\Omega.
\end{multline*}
\end{proof}

From the formulation of the main conjecture in \cite{HK1} we can
deduce the following weaker instance of Theorem~\ref{thm:Main
Theorem}.(ii).

\begin{lem}\label{lem:main normalisiert}
Let $\phi:\Omega\To\widetilde{\Omega}$ be the normalisation of
$\Omega$. Then
$$
\CHAR\LeftD\phi_*\Rc(U_{K_\infty},\ValR(\rho))=L(U_{K_\infty},\rho^{-1}\cyclchar)\widetilde{\Omega}.
$$
\end{lem}
\begin{proof}
We may decompose by the characters of $\GAL(K_0/\Rat)$. Thus, we
may assume that $K_\infty=\Rat_\infty$. By Lemma~\ref{lem:indep of
Euler-f} we can further reduce to $U=\SPEC \algebra{\Int}{1/p}$.

If $\rho$ is odd, then $L(U_{\Rat_\infty},\rho^{-1}\cyclchar)=1$
by definition. On the other hand,
$$\CHAR \Rc(U_{\Rat_\infty},\ValR(\rho))=(1)$$
by \cite{HK1}, Theorem~4.2.2. If $\rho$ is even, then
$c(U_{\Rat_\infty},\rho)=0$ and by \cite{HK1}, Theorem~4.2.4,
$$
\CHAR\RightD\IwaSect(U_{\Rat_\infty},\ValR(\rho))
=L(U_{\Rat_\infty},\rho^{-1}\cyclchar)\power{\ValR}{\GAL(\Rat_{\infty}/\Rat)}.
$$
Originally, both theorems only deal with the case that $\rho$ is a
finite character times an integral power of $\cyclchar$, but the
general case follows easily by twisting.
\end{proof}

We will now turn to the proof of Theorem~\ref{thm:Main
Theorem}.(i). A large portion of it can be dealt with by the
following lemma. Here, we use the Theorem of Ferrero-Washington
for the second time (see Proposition~\ref{prop:mu=0 for
L-elements}).

\begin{lem}\label{lem:even part}
$\Tw{\rho^{-1}\cyclchar}(\pr{-})\IwaH^1(U_{K_\infty},\ValR(\rho))$
and $\IwaH^2(U_{K_\infty},\ValR(\rho))$ are finitely generated as
$\ValR$-modules.
\end{lem}
\begin{proof}
By Corollary \ref{cor:modified_compatibility} we may enlarge
$K_\infty$ such that $\rho$ factors through $\GAL(K_\infty/\Rat)$.
Further, nothing changes if we then twist by $\rho^{-1}\cyclchar$.
Hence, we may assume $\rho=\cyclchar$.

Set $X=\SPEC\Int$, $S=X-U$. By the Theorem of Ferrero-Washington
(\cite{Wa1}, Theorem~7.15) the module
$$
\varprojlim_n \ValR\tensorp{\Int}\Pic(X_{K_n})
$$
is finitely generated over $\ValR$. By Proposition~\ref{prop:gysin
factors} this is also true for the modules
$$
\varprojlim_n \etH^i(S_{K_n},\ValR).
$$
Further, it is an elementary fact of the theory of cyclotomic
fields that
$$
\card{(\Gm(X_{K_n})/\mu(X_{K_n})\pr{+}\Gm(X_{K_n}))}\leq 2,
$$
where $\mu$ denotes the sheaf of unit roots (see \cite{Wa1},
Theorem~4.12). In particular, as $p$ was assumed to be an odd
prime,
$$
\varprojlim_n \ValR\tensorp{\Int}\pr{-}\Gm(X_{K_n})= \varprojlim_n
\ValR\tensorp{\Int}\mu(X_{K_n}).
$$
This module is obviously finitely generated over $\ValR$ as well.
Now use the exact sequence of Proposition~\ref{prop:cohom of
Zp(1)}.(ii).
\end{proof}

After decomposition by characters we may assume that $K_\infty$ is
a $p$-extension; in particular, totally real. Additionally, we may
shrink $U$ by Lemma~\ref{lem:indep of Euler-f} such that
$U_{K_\infty}$ and $\rho$ are unramified over $U$. The case that
$\rho$ is even has already been settled by the above lemma. The
key to the remaining case is the following

\begin{lem}\label{lem:key lemma}
Let $K_\infty$ be a $p$-extension, $\rho$ be odd, and let both be
unramified over $U$. Write $\ideal{p}$ for the prime ideal of
$\Omega$ with $p\in\ideal{p}$ and $\CODIM \ideal{p}=1$. Then there
exists a nonzero divisor $x$ of $\Omega_\ideal{p}$ and a
quasi-isomorphism
$$
\Rc(U_{K_\infty},\ValR(\rho))_{\ideal{p}} \isomorph
\Omega_{\ideal{p}}/x\Omega_{\ideal{p}}[-1].
$$
\end{lem}
\begin{proof}
We will first show that
$\RightD\IwaSect(U_{K_\infty},\ValR(\rho))$ is quasi-isomorphic to
a complex $\cmplx{P}$ of finitely generated projective
$\Omega$-modules with $P^i=0$ for $i\notin\set{1,2}$. By
Proposition~\ref{prop:cohom dim} we can achieve that $P^i=0$ for
$i>2$.

Recall that in the present situation, $\Omega$ is a local ring.
Let $k$ be the residue field of $\Omega$.
Proposition~\ref{prop:tensor products} implies that
$$
k\dtensorp{\Omega}\RightD\IwaSect(U_{K_\infty},\ValR(\rho))=
\RightD\etSect(U,k\tensorp{\Omega}\Ind\ValR(\rho)).
$$
Since $K_\infty$ is totally real and $\rho$ is odd, every lift of
the complex conjugation will act by multiplication by $-1$ on
$\Ind\ValR(\rho)\tensorp{\Omega}k$; consequently,
$$
\etH^0(U,k\tensorp{\Omega}\Ind\ValR(\rho))=0.
$$
Thus, we can choose $P^i=0$ for $i<1$ as well.

Lemma~\ref{lem:even part} then implies that the complex
$\RightD\IwaSect(U_{K_\infty},\ValR(\rho))_{\ideal{p}}$ is
quasi-isomorphic to a free $\Omega_{\ideal{p}}$-module sitting in
degree $1$. The claim follows since $\Omega c(U_{K_\infty},\rho)$
is free of rank $1$  and $\Rc(U_{K_\infty},\ValR(\rho))$ is
torsion by Lemma~\ref{lem:torsion+free}.
\end{proof}

Putting this and Lemma~\ref{lem:main normalisiert} together we see
that in the situation of Lemma~\ref{lem:key lemma}, the invertible
ideals of the normalisation of $\Omega_\ideal{p}$ generated by
$x^{-1}$, respectively $L(U_{K_\infty},\rho^{-1}\cyclchar)$,
agree. But
$$
L(U_{K_\infty},\rho^{-1}\cyclchar)=1;
$$
hence, $x$ is unit in the normalisation of $\Omega_{\ideal{p}}$
and therefore a unit in $\Omega_{\ideal{p}}$ itself. This finishes
the proof of Theorem~\ref{thm:Main Theorem}.(i).

At last, we complete the proof of Theorem~\ref{thm:Main
Theorem}.(ii). Let $(K_\infty,\rho,U)$ be any admissible triple.
By Proposition~\ref{prop:comparison of ideals} it suffices to show
that
$$
\CHAR \Rc(U_{K_\infty},\ValR(\rho))_{\ideal{p}}=
L(U_{K_\infty},\rho^{-1}\cyclchar)\Omega_{\ideal{p}}
$$
for all prime ideals $\ideal{p}$ of codimension $1$. In
Lemma~\ref{lem:main normalisiert} we have already proved this for
those $\ideal{p}$ that do not contain $p$. By
Theorem~\ref{thm:Main Theorem}.(i) and Proposition~\ref{prop:mu=0
for L-elements} the equality also holds for the remaining primes.
\end{proof}

\bibliographystyle{amsalpha}
\bibliography{IwaBib}
\end{document}